\pretolerance=2000
\hbadness=2000
\vbadness=2000

\parindent=1truecm

\def\id{1\kern-2.5pt{\rm l}}

\def\compl{{\rm l\!\!\!C}}

\def\tr{{\rm Tr}}

\def\gl{\hbox{\got gl}}

\def\so{\hbox{\got so}}

\font\got=eufm10
\font\eusm=eusm10
\font\ccap=cmbx12 at 18pt
\font\scap=cmbx12 at 12pt

\font\sscap=cmbxsl10

\def\S{\hbox{\got S}}
\def\A{\hbox{\eusm A}}
\def\H{\hbox{\got H}}
\def\M{\hbox{\got M}}

\def\g{\hbox{\got g}}
\def\G{\hbox{\eusm G}}
\def\ad{\rm ad}

\def\gl{\hbox{\got gl}}

\vskip 2pt plus 0.5fill
\centerline{\ccap On a Poisson structure}
\medskip
\centerline{\ccap on the space of Stokes matrices }
\bigskip
\vskip 2pt plus 0.2fill
\centerline{\scap Monica Ugaglia}\medskip
\centerline{\sscap Scuola Internazionale Superiore di Studi Avanzati}
\centerline{\sscap Via Beirut 4, 34014 Trieste, Italy}
\vskip 4truecm
\noindent{\bf Abstract:}
\noindent In this paper we study the map associating to a linear differential 
operator with rational coefficients its monodromy data. The operator is
of the form  $\Lambda(z)={d\over dz}-U-{V\over z}$, with one regular and 
one irregular singularity of Poincar\'e rank 1, where $U$ is a diagonal and 
$V$ 
is a skewsymmetric $n\times n$ matrix. We compute the Poisson structure of
the corresponding
 Monodromy 
Preserving Deformation Equations (MPDE) on the space of the monodromy data.
\vskip 9truecm
\centerline{Preprint SISSA 120/98/FM}
\vskip 5truecm
\vskip 0pt plus 0.5fill
\vfill\eject

\vskip 2truecm
\centerline
{\scap  0. Introduction}
\bigskip
\bigskip
\noindent Monodromy preserving deformation equations (MPDE) of linear differential
operators with rational coefficients are known since the beginning of the
century [Fu, Schl, G]. Particularly, the famous six
Painlev\'e equations are known [G] to be of this type. MPDE were included
in the framework of the general theory of integrable systems much later,
at the end of 70s [ARS, FN1, JMU]; see
also [IN]). Many authors were inspired by the parallelism
between the technique of soliton theory based on the spectral transform
and that of the MPDE theory based on the monodromy transform.
Another issue of this parallelism between soliton equations and MPDE
is that, in both cases, one deals with certain classes of Hamiltonian
systems, namely, with infinite-dimensional Hamiltonian structures of
evolutionary equations and of their finite-dimensional invariant
submanifolds in soliton theory, and with remarkable finite-dimensional
time-dependent Hamiltonian systems in the MPDE theory.

Recall that one of the first steps in soliton theory was understanding
of the Hamiltonian nature of the spectral transform as the transformation
of the Hamiltonian system to the action-angle variables [ZF].
Further development of these ideas was very important for development
of the Hamiltonian approach to the theory of solitons [FT]
and for the  creation of a quantum version of this theory.

In the general theory of MPDE it remains essentially an open question
to understand the Hamiltonian nature of the monodromy transform, i.e., of
the map
associating the monodromy data to the linear differential operator with
rational coefficients. This question was formulated in [FN2]
and solved in an example of a MPDE of a particular second order
linear differential operator. However, the general algebraic properties
of the arising class of Poisson brackets on the spaces of monodromy
data remained unclear. The technique of [FN2] seems not to work
for more general case. The authors of the papers [AM, FR, KS, Hi]
consider the important case of MPDE of Fuchsian systems in a more general
setting of symplectic structures on the moduli space of flat connections
(see, e.g., [A]) not writing, however, the Poisson bracket on the
space of monodromy data in a closed form. MPDE of non-Fuchsian operators
and Poisson structure on their monodromy data were not considered in these
papers.

In the present paper we  solve the problem of computing  the Poisson
structure of MPDE in the monodromy data coordinates for one particular
example of the operators with one regular and one irregular singularity
of Poincar\'e rank 1
$$
\Lambda(z)={d\over
dz}-
U-{V\over z}
$$ 
where $U$ is a diagonal matrix with pairwise distinct entries and
$V$ is a skewsymmetric $n\times n$ matrix.
Recently MPDE of this operators proved to play a fundamental role
in the theory of Frobenius manifolds [D, D1]. The Poisson structure
of MPDE for the operator $\Lambda$ coincides with the standard
linear Poisson bracket on the Lie algebra $\so(n)\ni V$. The most
important part of the monodromy data is the Stokes matrix (see
the definition below). This is an upper triangular matrix $S=S(V,U)$ with
all
diagonal entries being equal to 1. Generically $S$ determines other
parts of the monodromy data. It turns out that, although the monodromy map
$$
V\mapsto S
$$
is given by complicated transcendental functions,  
the Poisson bracket on the space of Stokes matrices is given by
very simple degree two polynomials (see formula (3.2) below).
 The resulting Poisson bracket does not depend on $U$ since is involved in 
the Hamiltonian description of the isomonodromy deformations of the operator $
\Lambda(z)$. The
technique of [KS] was important in the derivation of this main
result of the present paper.

We hope that this interesting new class of polynomial Poisson brackets
and their quantization (cf. [R, Ha2]) deserves a further
investigation that we are going to continue in subsequent publications.

The paper is organized as follows: after  recalling some basic
notations,
in section 1.1 we describe the monodromy of the operator
$\Lambda(z)$ 
around the two singular points; in section 1.2 we present the
MPDE for this operator. In section 2.1 and 2.2  we describe 
the related Fuchsian 
system and its MPDE; in section 2.3 the Poisson structure on the space of
monodromy data of the Fuchsian system is described.
In section 3 we give the relation between the monodromy data of
the two systems and we explicitly calculate the Poisson bracket
on the space of the Stokes matrices.
\bigskip
\vfill\eject
\bigskip

\bigskip
\noindent{\sscap 0.1 Basic notations}
\bigskip

\noindent Let us consider in the complex domain a differential equation with
rational coefficients
$$
{{dy}\over{dz}}=A(z)y(z)\eqno{(0.1)}
$$
where $$y=\left(\matrix{y_1\cr y_2\cr ..\cr y_n\cr}\right),\qquad 
A(z)=\left(\matrix
{a_{11}&a_{12}&..&..&a_{1n}\cr
 a_{21}&..&..&..&a_{2n}\cr
 ..&..&..&..&..\cr
 a_{n1}&..&..&..&a_{nn}\cr}\right)$$

An arbitrary solution $y(z)$ of $(0.1)$ is locally holomorphic but
globally 
multivalued; the poles of $A(z)$ are 
singularities of the solution. Fixing a basis 
$y^{(1)},\ldots,y^{(n)}$ in the $n$-dimensional space of solutions we 
construct the {\it fundamental} $n\times n$ matrix

$$
Y(z)=(y^{(1)},\ldots,y^{(n)})
$$

\noindent satisfying the matrix version of (0.1)

$$
{{dY(z)}\over{dz}}= A(z) Y(z).\eqno{(0.2)}
$$
\bigskip

\centerline {\scap 1. Systems with irregular singularity}
\bigskip
\noindent{\sscap 1.1 Stokes phenomenon}
\bigskip
\noindent In this paper we will concentrate our attention on the 
linear systems
$$
{dY\over {dz}}=(U+{V\over{z}})Y,\qquad z\in\compl,\eqno{(1.1)}
$$
where $U$ is a diagonal $n\times n$
matrix with distinct entries $u_1,u_2,\ldots,u_n$ and $V=(v_{ij}) \in
\so(n,\compl)$, with  nonresonant eigenvalues $(\mu_1, \mu_2,\dots,\mu_n)$
(i.e.  $\mu_i-\mu_j\not\in Z\setminus 0$). The
solutions of the system (1.1) have
two singular points, $0$ and $\infty$.
 \bigskip
\noindent $\bullet$ Near the  point $z=0$  a  fundamental matrix of solutions
$Y_0(z)$ exists such that 
 $$Y_0(z)=W(z)z^{\theta}=[W_0+W_1
z,\ldots]z^{\theta},\eqno{(1.2)} 
$$
where $\theta$ is the diagonalization of $V$,\ \  $\theta =W_0^{-1}VW_0
= diag (\mu_1, \mu_2, \dots, \mu_n),$ 
and  $W(z)$ converges for small $|z|$.
Such  kind of singularities is called {\it Fuchsian}.

If one continues $Y_0(z)$ along a path encircling 
the point $z=0$, the columns of the resulting matrix are linear
combinations of the columns of $Y_0(z)$; hence there exists a matrix
$M_0$ such
that
$$
Y_0(z)\mapsto Y_0(z)M_0.
$$
 The matrix $M_0$ is called {\it monodromy} matrix around
zero. In our case  $M_0=exp (2\pi i
\theta)$.  
\bigskip 
\noindent $\bullet$ At $\infty$ the solution has an  {\it irregular} singularity
of  Poincar\'e rank 1. This means that
 it is possible to construct a formal
series 
$$
\Gamma(z)=\id
+{\Gamma_1\over{z}}+{\Gamma_2\over{z^2}}+\dots
$$
where $V=[\Gamma_1,U]$+ {\it diagonal}, i.e.
$\Gamma_1=(\gamma_{ij})=({v_{ij}\over{u_j-u_i}})$ for $i\neq j$,
 and to define certain sectors $\S_i$ in which  a fundamental
matrix of solutions $Y_i$ exists with asymptotic behavior 
$$
Y_i\sim Y_{\infty}=\Gamma(z)e^{zU}, \eqno{(1.3)} 
$$
for $\vert z\vert\rightarrow\infty$ in $\S_i$.
This means that $\Gamma(z)$ is the asymptotic expansion of
$
Y_ie^{-zU}.
$

In different sectors one has different solutions, and this fact is known as
{\it Stokes phenomenon}.
The matrices connecting the solutions in different sectors are called
{\it Stokes
matrices}.
\bigskip
\noindent
A complete and detailed description of the phenomenon can be found in
[BJL1], [Si],
[IN],[U];  here we will concentrate our attention on the particular
operator $\Lambda(z)$ {\bf = ${d\over dz} -U-{V\over z}$} (see also [D]). 

\bigskip
\noindent Following [D] we define 
an {\it admissible} line for the system (1.1) as a line $l$ 
through the origin on the $z$-plane such that
$$
{\rm Re} z(u_{_{i}}-u_{_{j}})\vert_{z\in l}\not=0\qquad
\forall i\not=j.
$$
We denote the half-lines
$$
l_+={z:\arg z=\psi}\qquad l_-={z:\arg z=\psi-\pi},
$$
where $\psi$ is a fixed real value of the argument.

\bigskip
\noindent The line $l$ lies in the intersection $\S_+\cup \S_-$ of the two sectors
$$\S_R:\ \psi-\pi-\epsilon< \arg z<\psi+\epsilon$$
 and $$\S_L:\ 
\psi-\epsilon< \arg z<\psi+\pi+\epsilon.$$
Here $\epsilon$ is a sufficiently small positive number.
\bigskip

\noindent{\bf Theorem 1.1} : {\it There exists a unique  solution
$Y_{L}(z)$  analytic in the sector $\S_L$  with the asymptotic
behavior} $$
Y_{L}(z)\sim Y_{\infty};  
$$
{\it the same holds for $Y_{R}(z)$ in $\S_L$}.
\bigskip
\noindent{\bf Proof}: See  [BJL1].
\bigskip 
\noindent   $S_+$ and $S_-$ are the  Stokes matrices connecting the two
solutions in $\S_+$, resp.  in $\S_-$, i.e.
$$
Y_L(z)=Y_R(z)S_+, \qquad z\in \S_+
$$
and
$$
Y_L(z)=Y_R(z)S_-, \qquad z\in \S_-.
$$
 From the skew-symmetry $V^T=-V$ it follows
$$
S_-=S_+^T.
$$ 

Moreover, one can prove that,  given an admissible line, it is possible
to order the elements $u_i$, i.e., to perform a conjugation
$$
\Lambda(z)\mapsto P^{-1}\Lambda(z) P,
$$
where $P$ is the matrix of the permutation
 in such a way that 
the  Stokes matrix $S\equiv S_+$ is upper triangular.
\bigskip 
\noindent{\bf Remark}:  The full set of monodromy data for the operator
$\Lambda$ consists of
the Stokes matrix $S$ but also of the monodromy
matrix at the point $0$ and of the matrix $C$ connecting the
solution (1.2) near
zero with a solution near the infinity:
$$
Y_0(z)=Y_L(z) C.
$$
The
monodromy data $\{S,M_0,C\}$ satisfy
certain constraints described in [D1]. Particularly, 
$$C^{-1} {S^T}^{-1} S \, C =M_0.
$$
So, in the generic case (i.e., the diagonalizable and nonresonant one)
under consideration the diagonal entries of $M_0$ $e^{2\pi i\mu_1}$, 
\dots, $e^{2\pi i \mu_n}$ are the eigenvalues of ${S^T}^{-1}S$ and $C$ is
the diagonalizing transformation for this matrix. The ambiguity in the
choice of the diagonalizing transformation does not affect the operator
$\Lambda$. So, the ${n(n-1)\over 2}$ entries of the Stokes matrix $S=S(V;U)$ 
can
serve as local coordinates near a generic point of the space of monodromy
data of the operator $\Lambda$ (see details in [D], [D1]).
\bigskip
\bigskip 
\noindent{\sscap 1.2 Monodromy Preserving Deformation Equations}
\bigskip
\noindent  
MPDE describe
 how should
the matrix $V$ be
deformed, as a function of the ``coordinates" $u_i$, in order to preserve the
monodromy data. MPDE are the analogue of the isospectral equations
in soliton
theory.
\noindent The MPDE for the operator
$\Lambda (z)={d\over{dz}}-U-{V\over{z}}$ are obtained (see [U], [D]) as
compatibility
equations of the system (1.1) with the system
$$
{\partial Y\over \partial u_i}= \bigl(zE_i-V_i\bigr) Y,
$$
 where
$V_i=[\Gamma_1,E_i]=ad_{_{E_i}}ad_{_{U}}^{-1}V$ and $(E_i)^a_b=\delta^a_i
\delta^i_b$. These equations admit the Lax form

 $$ {\partial V\over {\partial u_i}}=[V,V_{i}].\eqno{(1.4)}  
$$
In the generic case (see Remark above) the solution $V=V(U)$ of the MPDE 
can be locally written in implicit form
$$
S(V;U)=S\eqno{(1.5)}
$$
for a given constant Stokes matrix $S$. In other words, the entries of the 
Stokes matrix can serve as a complete system of first integrals of the MPDE 
(1.4). To explicitly resolve the system (1.5) one has to solve an appropriate 
Riemann--Hilbert boundary value problem. Although this can be  
 explicitly done in a very few cases, one can extract certain important
information regarding the analytic properties of the solution;
see more detailed discussions of these properties in [IN], [JM], [JMU], [Si].

One can write the MPDE as  a Hamiltonian system on the
space of the
skewsymmetric matrices $V$ with the standard linear Poisson bracket
for $V=(v_{ab})\in\so (n)$: 
$$
\{v_{ab},v_{cd}\}=v_{ad} \delta_{bc}+v_{bc} \delta_{ad}-v_{bd} \delta_{ac}-
v_{ac} \delta_{bd}.\eqno{(1.6)}
$$
Indeed, the Lax equation (1.4) can be rewritten as 
$$
{\partial V\over {\partial u_i}}=\{V,H_i(V,u)\},
$$
for the Hamiltonian function 
$$H_i={1\over 2}\sum_{j\not=i}{v^2_{ij}\over
u_i-u_j}.\eqno{(1.7)}
$$
\bigskip
\noindent In this case, the Poisson bracket is linear but the dynamic of the
problem is very complicated; in the following we will show how, very much as in
the case of isospectral equations, it is possible to find a different coordinate
system (the entries of the Stokes matrix) in which the dynamic of the evolution
is trivial, but the Poisson structure is quadratic. 
The technique developed here
consists in building up the monodromy map $V\to S$ passing through an
auxiliary 
Fuchsian system. 
The MPDE for the system (1.1) can be represented also as MPDE for an
appropriate
Fuchsian system  
$$
{d\chi\over d\lambda}=\sum_{i=1}^n{A_i\over \lambda-u_{_{i}}}\chi,
$$
which we shall describe in the next
section. The basic idea to construct the Poisson bracket on the space of 
Stokes matrices is to include the  map from $V\in \so(n)$ to $S\in{\cal S}$ 
into the
following commutative
diagram of Poisson maps
$$
\eqalignno{&\so(n)\longrightarrow\ {\cal S}\cr
 &\downarrow\qquad\qquad\downarrow&(1.8)\cr
&\A/\G \longrightarrow
\M/GL(n,\compl)\cr}
$$
where $\A/\G$ is the space of residues $\{A_i\}$ of the connection 
$A=\sum_{i=1}^n{A_i\over \lambda-u_{_{i}}}d\lambda$ modulo the action of
the gauge
group $\G$, as we will explain in section 2.2, and $\M/GL(n,\compl)$ is
the space
of the monodromy data of the
Fuchsian system
(section 2.3), i.e. the space of $n$-dimensional representations of the
free group with $n$ generators.
\bigskip
\bigskip
\centerline{\scap 2. Related Fuchsian system}
\bigskip
\noindent{\sscap 2.1 Fuchsian system}
\bigskip
\noindent One can relate the system (1.1), with one regular  and one irregular
singularity to a
 system with $n+1$ Fuchsian singularities:
$$
{d\Phi\over d\lambda}=\sum_{i=1}^n{B_i\over \lambda-u_{_{i}}}\Phi,\eqno{(2.1)}
$$
where
$$
B_i=-E_i\biggl(V+{1\over2}\id\biggr),\quad i=1,\ldots,n 
$$
and
$$
B_{\infty}=V+{1\over2}\id.
$$
Such a relation is well known in the domain of differential equations, see, e.g.
[BJL], [Sch].
\bigskip
\noindent Now we will briefly describe the monodromy data of the system (2.1).
\bigskip
\noindent In this case $u_j$ is a Fuchsian singular points and, as in (1.2),
the general solution near $u_j$ can be expressed as 
$$
\Phi_j(\lambda)=W^{^{(j)}}(\lambda) (\lambda-u_j)^{^{\hat
B_j}},
$$
where $W^{^{(j)}}(\lambda)=W^{^{(j)}}_0+(\lambda-u_j)W^{^{(j)}}_1 +\ldots$
 converges  for small $|\lambda-u_j|$ and 
$\hat B_j = 
-{1\over 2} E_j$ is the
diagonalization of $B_j$.

\noindent We  denote  $M_j$
the monodromy matrix along the path $\gamma_j$ 
encircling
 the point $u_j$ w.r.t. the basis $\Phi_\infty$ we define in (2.2) below.
The matrix $M_j$ is conjugated with the matrix $exp (2\pi i\hat B_j)$.

\noindent Also the point $\infty$ is Fuchsian; the general solution can be 
expressed as
$$
\Phi_{\infty}(\lambda)=W^{^{(\infty)}}(\lambda)
\bigl({1\over\lambda}\bigr)^{^{\hat B_{\infty}}},
\eqno{ (2.2)} $$
where $W^{^{(\infty)}}(\lambda)=W^{^{(\infty)}}_0+{W^{^{(\infty)}}_1
\over\lambda} +\ldots$ converges at $|\lambda| \to \infty$ and
$\hat B_{\infty}= diag( {1\over 2}+\mu_1, \dots, {1\over 2}+\mu_n) $ is
the diagonalization of $B_{\infty}$. Indeed, the following
relation holds in the space of the residues :
$$
-\sum_{i=1}^n B_i=B_{\infty} = {1\over 2}\id+V.
$$
In this basis the monodromy matrix $M_\infty =-e^{2\pi i \theta}$. We
assume that the loops $\gamma_1$, \dots, $\gamma_n$ and $\gamma_\infty$
are chosen in such a way that
$$
M_1 M_2\dots M_n M_\infty=1.\eqno{(2.3)}
$$

\bigskip 
\noindent{\sscap 2.2 Monodromy Preserving Deformation equations}
\bigskip
\noindent We now  want to deduce the MPDE
 for the system $(2.1)$. This amounts to find
how can the matrix $B_j$ be deformed as function of $u_1,u_2,\ldots,u_n$ in 
order to preserve the monodromy matrices $M_1$, \dots, $M_n$, $M_\infty$.
The answer is given by 
 
\noindent{\bf Theorem 2.1 (Schlesinger)}: {\it If the fundamental solution
 near 
 infinity is normalized as in (2.2) and $A_{\infty}$ is a constant 
diagonal matrix with 
nonresonant elements,  then the
 dependence 
of the $A_j$ on 
the position of the poles of the Fuchsian system} 

$$
{d\Phi\over d\lambda}=\sum_{i=1}^{n}{A_i\over \lambda-u_{_{i}}}\Phi
$$
{\it is given, in order to preserve the monodromy,  by} 
$$
\cases{{{\partial A_i}\over{\partial u_j}}={{1}\over{u_i-u_j}}[A_i,A_j]\ 
\ \ i\not=j&\cr {{\partial A_j}\over{\partial u_j}}=-\sum_{i\not= 
j}{{[A_i,A_j]}\over{u_i-u_j}},&\cr}
$$
{\bf Proof}: it can be found in  [Si].
\bigskip
Note that system (2.1) does not satisfy the hypotheses of the Schlesinger
theorem, because
$B_{\infty}=\bigl(V+{1\over 2}\id\bigr)$ is not diagonal. 

In order to apply the
Schlesinger theorem  it is sufficient to  perform the gauge transformation
$$
B_i\mapsto A_i=W_0^{-1}B_iW_0,\eqno{(2.4)}
$$
where $W_0$ is the matrix of eigenvectors of $V$ 
normalized in such a 
way that
$$
{\partial W_0\over\partial u_i}=\ad_{E_i}\ad_U^{-1}V.\eqno{(2.5)}
$$
Indeed,  substituting $\Phi=W_0\chi$, the system (2.1) transforms into
$$
{d\chi\over d\lambda}=\sum_{i=1}^n{A_i\over \lambda-u_i}\chi\eqno{(2.6)}
$$
 and the Schlesinger system follows from the compatibility of (2.6) 
 with
$$
{\partial\chi\over\partial u_i}= -{A_i\over \lambda-u_i}\chi.
$$
(See [D]).
\bigskip

\noindent The Schlesinger system  can be rewritten in the Hamiltonian form 
$$
{d A_i\over d u_j}=\{A_i,\H_j\}
$$
with  the Hamiltonians
$$
\H_j=-\sum_{k\not=j}{\tr(A_jA_k)\over u_j-u_k}
$$
w.r.t.  the linear Poisson bracket
$$
\{(A_i)^a_b,(A_j)^c_d\}= \delta_{ij}\bigl(
\delta^a_d(A_i)^c_b-\delta^c_b(A_j)^a_d\bigr).\eqno{(2.7)}
$$
This corresponds to taking, for every $u_i$, the residue
$A_i\in\gl(n,\compl)$
with the natural Poisson bracket on  $\gl(n,\compl)$. The residues relative to
different singular points commute. In other words (see [KS],[FR],[A]) this
corresponds  to  read the matrices $A_i$ as residues of a
flat connection (with values in the Lie algebra $\g=\gl(n,\compl)$) on the
Riemann
surface with $n+1$ punctures: 
 $$
A=\sum_{i=0}^n{A_i\over\lambda-u_i}d\lambda
$$
(in our case $u_0=\infty$).
On the space of flat connections
modulo gauge
transformations it is defined the Poisson bracket
$$
\{A^{^{a}}(\mu),A^{^{b}}(\nu)\}=-f^{^{ab}}_{_{c}}{A^{^{c}}(\mu)-
A^{^{c}}(\nu)\over \mu-\nu},
$$
where
 $f^{^{ab}}_{_{c}}$ are the structure constants of $\g$ w.r.t. the basis 
$\{e_{_{a}}\}$ and
$$
A^{^{a}}(\mu)=\sum_{i=0}^n{A^{^{a}}_i\over\mu-u_i},\qquad A(\mu)
=A^{^{a}}(\mu) e_{_{a}}.
 $$
 
 This
Poisson bracket gives (2.6).
\bigskip
\noindent Now we can perform the first step in the construction of the map between $V$ and
$S$, that is we pass from  $\so(n)$  to the space $\A/\G$, where
$$
\A=\bigl\{V,A_1,\ldots ,A_n \vert \sum_{i=0}^n A_i=0\bigr\}
$$
is the family of the residues of $A(\lambda)$ and $\G$ is the gauge
group.
\bigskip
\noindent {\bf Lemma 2.1 :}{\it The map
$V\in\so(n)\mapsto (V,A_1,\ldots,A_n)\in\A/\G$ is a Poisson map.} (Cf. [Ha1],[Hi]) 
\bigskip

\noindent {\bf Proof}: We must compare the Poisson brackets on the two spaces. In
$\so(n)$ 
one has the natural coordinates $\{v_{ab}\}$, with the Poisson bracket
(1.6).
The natural coordinates  in the quotient space
$\A/\G$ are
the traces of the products of 
the  matrices $A_i$, so that we consider
the brackets
$$
\eqalignno{\{\tr (A_i A_k)&,\tr (A_j A_l)\}=\{(A_i)^a_b(A_k)^b_a,
(A_j)^c_d(A_l)^d_c\}=\cr
&=(A_i)^a_b(A_j)^c_d\{(A_k)^b_a,(A_l)^d_c\}+
(A_i)^a_b(A_l)^d_c\{(A_k)^b_a,(A_j)^c_d\}+\cr
&+(A_k)^b_a(A_j)^c_d\{(A_i)^a_b,(A_l)^d_c\}+
(A_k)^b_a(A_l)^d_c\{(A_i)^a_b,(A_j)^c_d\}&(2.8a)
 \cr} 
$$
and
$$
\eqalignno{\{\tr (A_i V),&\tr (A_j V)\}=\{(A_i)^a_bV^b_a,
(A_j)^c_dV^d_c\}=\cr
&=(A_i)^a_b(A_j)^c_d\{V^b_a,V^d_c\}+
V^b_aV^d_c\{(A_i)^a_b,(A_j)^c_d\}.&(2.8b)
 \cr} 
$$
On $\A/\G$ by direct calculation, using the bracket (2.7), one obtains
$$
\eqalignno{\{\tr (A_i A_k),\tr (A_j A_l)\}&=\delta_{kl}\tr\bigl(
A_iA_jA_k-A_kA_jA_i\bigr)+\delta_{kj}\tr\bigl(
A_iA_lA_k-A_kA_lA_i\bigr)+\cr
&+\delta_{il}\tr\bigl(
A_kA_jA_i-A_iA_jA_k\bigr)+\delta_{ij}\tr\bigl(
A_kA_lA_i-A_iA_lA_k\bigr)=\cr
&=2\biggl(\delta_{kl}\tr\bigl(
A_iA_jA_k\bigr)+\delta_{kj}\tr\bigl(
A_iA_lA_k\bigr)
+\delta_{il}\tr\bigl(
A_kA_jA_i\bigr)+\delta_{ij}\tr\bigl(
A_kA_lA_i\bigr)\biggr).&(2.9)}
$$
Indeed, $A_i=-E_i\bigl(V+{1\over2}\id\bigr)$ implies 
$$
\tr\bigl(A_iA_jA_k\bigr)=-\tr\bigl(A_kA_jA_i\bigr)=v_{ij}v_{jk}v_{ki}.
\eqno{(2.10)}
$$

On the other hand, $\tr(A_iA_j)=-v_{ij}^2$, hence
$$
\eqalignno{\{\tr (A_i A_k),\tr (A_j A_l)\}&=4v_{ik}v_{jl}\{v_{ik},v_{jl}\}\cr
&=4\bigl(\delta_{kl} v_{ij}v_{jk}v_{ki}-\delta_{kj} v_{ik}v_{kl}v_{li}+
\delta_{il} v_{ik}v_{kj}v_{ji}+\delta_{ij} v_{ik}v_{kl}v_{li}\bigr)&(2.11)}
$$
where we have used the bracket (1.6). By means of (2.10) it is easy to 
check that it coincides with (2.9).
\bigskip

The same can be done for equation (2.8b). Indeed, using the bracket on 
the $A_i$ matrices and observing that 
$$
\tr\bigl(A_iA_jV\bigl)=-\tr\bigl(VA_jA_i\bigr)=\sum_{k\not=i\not=j}
v_{ij}v_{ki}v_{jk},
$$
one finds
$$ 
\{\tr (A_i V),\tr (A_j V)\}=2\tr (A_iA_jV).
$$
On the other hand $\tr(A_iV)=-\sum_{k\not=i}v_{ki}^2$, that gives
$$
\{\tr (A_i V),\tr (A_j V)\}=4\sum_{k\not=i}\sum_{l\not=j}v_{ki}v_{lj}
\{v_{ki},v_{lj}\}=-4 \sum_{k\not=i}v_{ki}v_{kj}v_{ij}
$$
which coincides with (2.8b).
\hfill{Q.E.D}
\bigskip
\bigskip
\noindent {\bf
Lemma 2.2 :} {\it The MPDE for the system (1.1) and its related Fuchsian system
coincide.}
\bigskip
\noindent{\bf Proof}: It follows immediately from Lemma 2.1 by a
straightforward
calculation  using (2.5), that MPDE for the Fuchsian system (2.6) 
after the gauge
transformation (2.4) coincide with (1.4). Actually, 
one can see that the pull back of the Hamiltonian
$$
\H_j=-\sum_{k\not=j}{\tr(A_jA_k)\over u_j-u_k} =-\sum_{k\not=j}
{\tr (B_j B_k)\over u_j-u_k}
$$
is exactely equal to $H_j$, as defined in (1.7).
\bigskip
\bigskip
\noindent{\sscap 2.3 Poisson structure on monodromy data}
\bigskip
\bigskip
In this section we will perform the second step  of our
construction, that is we
will map the Poisson structure of $\A/\G$ into the space of monodromy data of the
Fuchsian system; this is shown in the following well-known (see, e.g.,
[Hi])
\bigskip
\noindent{\bf Theorem 2.2:} {\it The monodromy
map
$$
\A/\G\rightarrow\M/SL(n,\compl)
$$
where $\M=\bigl\{M_0,M_1,\ldots,M_n\vert M_1M_2\ldots M_nM_0=\id\bigr\}$,
is a Poisson map.}
\bigskip
\noindent To actually compute the Poisson bracket on the space of 
monodromy data, i.e., on the space of $n$-dimensional representations of
the free group with $n$ generators we will use,
following [KS] (Th. 4.2), the following technique. We construct the
skewsymmetric
bracket 
$$
\eqalignno{\biggl\{(M_i)^a_b,(M_j)^c_d\biggr\}&=i\pi
\biggl((M_j M_i)^c_b\delta^a_d+(M_i M_j)^a_d\delta^c_b- 
(M_i)^c_b (M_j)^a_d-(M_j)^c_b (M_i)^a_d\biggr)\qquad i<j&(2.12a)\cr
\biggl\{(M_i)^a_b,(M_i)^c_d\biggr\}&=i\pi\biggl(
(M_i^2)^c_b\delta^a_d-(M_i^2)^a_d\delta^c_b\biggr).&(2.12b)\cr}
$$
 on the space $\M$ of the monodromy
matrices. As it was proved in [KS],
when restricted
to the space of representations $\M/SL(n,\compl)$, this bracket
defines a Poisson structure on the quotient induced by the monodromy
map. Observe that the eigenvalues of the matrices $M_i$ are the Casimirs
of the Poisson bracket, i.e., the functions Poisson commuting with all
others (see [KS]).
\vfill\eject
\bigskip
\bigskip
\centerline{\scap 3. Poisson structure on the Stokes matrices}
\bigskip

\bigskip
\noindent{\sscap 3.1 Connecting the monodromy data of the two systems}
\bigskip
\noindent In the previous section we have seen that the space of
 monodromy data of a Fuchsian system  carries a natural Poisson
stucture. In this section we will show that this structure induces a
Poisson bracket on the space of  Stokes matrices of the related system we
studied in chapter 1. To
this end we consider  the relation between the  monodromy matrices 
$M_1,M_2,\ldots,M_n$ of the Fuchsian system and the Stokes matrix $S$. 

In section 2.1 we claimed that the two systems
$$
{dY\over {dz}}=(U+{V\over{z}})Y
$$
and
$$
{d\Phi\over d\lambda}=\sum_{i=1}^n{A_i\over \lambda-u_{_{i}}}\Phi
$$
are related, in the sense that, (see Lemma 2.3), the MPDE for the
operator $\Lambda(z)={d\over dz}-U-{V\over{z}}$ can be represented also 
as MPDE for the
operator $A(\lambda)=
{d\over d\lambda}-\sum_{i=1}^n{A_i\over \lambda-u_{_{i}}}$.
 \bigskip 
 For a detailed analysis of the transform  connecting the
two system see [D1]; here we will concentrate our attention on the relation
between  the monodromy data
of the two systems.
  
\noindent Following Theorem 2.2, we are interested
in the quotient of the space of the monodromy data of the Fuchsian system
w.r.t. the $GL(n,\compl)$ conjugations. 
So, we can choose a particular basis of solution of the
system and work with the corresponding monodromy matrices.
\bigskip

\noindent{\bf Theorem 3.1 :} {\it Suppose that $\bigl(S+S^T\bigr)$ is
nondegenerate; then  there exists a unique basis of solutions 
 (which depends on the particular choice of the branchcuts in the 
complex
$\lambda$--plane)
$\{\Phi^{^{(j)}}(\lambda)\}$ of the Fuchsian system  (2.1), such that 

- Near $u_i$ the solution has the behavior}
$$
\Phi^{^{(i)}}_{_{a}}\sim{1\over
\sqrt{u_i-\lambda}}\delta^{^{i}}_{_{a}}. 
$$

- {\it the monodromy matrices
are reflections, i.e., going  around the singularity $u_i$  the solutions
transform as}
 $$
\eqalignno{\Phi^{^{(i)}}&\rightarrow -\Phi^{^{(i)}}\cr
\Phi^{^{(j)}}&\rightarrow \Phi^{^{(j)}}-2g_{ij}\Phi^{^{(i)}}\cr}
$$
{\it where $G=(g_{ij})= {1\over 2}\left(S+S^T\right)$ is the Gram matrix
of the following
invariant bilinear
form w.r.t. the chosen basis}  
$$
g_{ij}=\biggl(\Phi^{^{(i)}},\Phi^{^{(j)}}\biggr):=\Phi^{^{(i)^T}}
\biggl(U-\lambda\biggr)\Phi^{^{(j)}}.
$$
{\it Invariance means that $g_{ij}$ does not depend on $\lambda$ neither
on
$u_1,\dots,u_n$}.
\bigskip
\noindent {\bf Proof}: See  [D1], Th.5.3.

\bigskip
\noindent {\bf Remark}:  $\Phi$ and $Y_L$ are related by the
Laplace transform  
$$ (Y_{L})^{^{(j)}}_{a}(z)={-\sqrt z\over
2\sqrt\pi}\int_{\gamma_{(j)}} \Phi^{^{(j)}}_{_{a}}(\lambda)
e^{\lambda z}d\lambda 
$$ where  ${\gamma_{(j)}}$ is a  fixed path  in the
$\lambda$--plane; analogously for
$Y_R$.
 
\bigskip

\noindent In the $\{\Phi^{^{(j)}}(\lambda)\}$ basis the $i$--th monodromy matrix $M_i$
has the form $$
M_i=\left(\matrix{1&0&\ldots&0&\ldots&0\cr
0&1&\ldots&0&\ldots&0\cr
\vdots&\vdots&\ldots&\vdots&\ldots&\vdots\cr
-2g_{1i}&-2g_{2i}&\ldots&-1&\ldots&-2g_{ni}\cr
\vdots&\vdots&\ldots&\vdots&\ldots&\vdots\cr
0&0&\ldots&0&\ldots&1\cr}\right),
$$
$$
2 g_{ij}=2 g_{ji} = s_{ij} ~~ {\rm for}~ i<j.
$$
This is a reflection w.r.t. the hyperplane normal to the vector 
$\Phi^{^{(i)}}$.
\bigskip

\noindent The Coxeter identity ([B]) gives 
$$
M_1 M_2 \ldots M_n= -S^{-1}S^T.
$$

\bigskip

\noindent{\bf Lemma 3.1:} {\it The following relations hold (all the
indices are pairwise distinct)} 
$$
\eqalignno{\tr (M_i M_j)&=n-4+s_{ij}^2&(3.1a)\cr
\tr(M_kM_iM_jM_i)&=n-4+(s_{kj}-s_{ij}s_{ik})^2&(3.1b)\cr
\tr(M_i M_j M_l M_k)&=n-8+s_{ij}^2+s_{ik}^2+s_{il}^2+
s_{jk}^2+s_{jl}^2+s_{kl}^2-s_{ij}s_{ik}s_{jk}+\cr
&-s_{ik}s_{il}s_{kl}-s_{jk}s_{jl}s_{kl}-
s_{ij}s_{il}s_{jl}+s_{ij}s_{il}s_{jk}s_{kl}.&(3.1c)\cr}
$$
\bigskip
\noindent{\bf Proof:} The fact that the $M_i$ are reflections and that 
$S+S^T=2G$ geometrically reads into
$$
-2\cos \alpha_{ij}=s_{ij}
$$
where $\alpha_{ij}$ is the angle between the two hyperplanes normal to $
\Phi^{^{(i)}}$ and $\Phi^{^{(j)}}$.

 On the other hand, the products $M_i
M_j$ are rotations by the angle $2\alpha_{ij}$ and this provides the relation
(3.1a), indeed 
$$
\tr (M_i M_j)=n-2+2\cos (2\alpha_{ij})=n-4+s_{ij}^2.
$$
\bigskip
\noindent To obtain  relation (3.1b) we observe that the product $M_iM_jM_i$
is still a reflection, w.r.t. the mirror normal to the vector
$M_i(\Phi^{^{(j)}})$. This means that the product $M_kM_iM_jM_i$
is  a rotation by the  angle $2\beta$, where 
$$
-2\cos
\beta=\biggl(M_i(\Phi^{^{(j)}}),\Phi^{^{(k)}}\biggr)=(\Phi^{^{(j)}}-s_{ij}
\Phi^{^{(i)}},\Phi^{^{(k)}})=s_{kj}-s_{ij}s_{ik} $$
so that
$\tr(M_kM_iM_jM_i)=n-2+2\cos(2\beta)=n-4+(s_{kj}-s_{ij}s_{ik})^2$. 
\bigskip
Finally, (3.1c) can be obtained   directly 
in the case of the $4\times 4$ reflection matrices $M_i$.
Indeed, for ordered indices $i,j,k,l$, the Coxeter identity gives 
$$
M_iM_jM_kM_l=-S_{ijkl}^{-1}S_{ijkl}^T,
$$
where
$$
S_{ijkl}=\left(\matrix{1&s_{ij}&s_{ik}&s_{il}\cr
0&1&s_{jk}&s_{jl}\cr
0&0&1&s_{kl}\cr
0&0&0&1\cr}\right).
$$
An easy calculation gives the result.
 
\noindent The same result  holds also 
 in dimension $n>4$. Indeed, one can observe that, for every $n$,
the product of four matrices acts non trivially 
in the $4$--dimensional subspace spanned by the vectors normal
to the mirrors  of the reflections $M_i,M_k,M_l,M_j$.
It is equal to the identity in the orthogonal complement
to the $4$-dimensional subspace.

\hfill{Q.E.D.}

 \bigskip

\noindent Combining all the above facts, we can conclude our construction proving
the following main
\bigskip
\noindent{\bf Theorem 3.2:} {\it 1. The following formulae
$$
\eqalignno{
\{s_{ik},s_{il}\}&={i\pi\over 2}(2s_{kl}-s_{ik}s_{il})\qquad i<k<l&(3.2a)\cr
\{s_{ik},s_{jk}\}&={i\pi\over 2}(2s_{ij}-s_{ik}s_{jk})\qquad i<j<k&(3.2b)\cr
\{s_{ik},s_{kl}\}&={i\pi\over 2}(s_{ik}s_{kl}-2s_{il})\qquad i<k<l&(3.2c)\cr
\{s_{ik},s_{jl}\}&=0\qquad i<k<j<l&(3.2d)\cr
\{s_{ik},s_{jl}\}&=0\qquad i<j<l<k&(3.2e)\cr
\{s_{ik},s_{jl}\}&=i\pi(s_{ij}s_{kl}-s_{il}s_{kj})\qquad i<j<k<l.&(3.2f)\cr}
$$
define a Poisson bracket on the space ${\cal S}$ of Stokes matrices.

2. The monodromy map
$$
\so(n)\to {\cal S}
$$
associating the Stokes matrix $S\in {\cal S}$ to the operator
$\Lambda ={d\over dz} -U - {V\over z}$, $V\in \so(n)$, is a Poisson map.

3. The eigenvalues of $S^{-1}S^T$ are the Casimir functions
of the Poisson bracket.

4. The Poisson bracket (3.2) is invariant w.r.t. the action of the braid
group $B_n$ on the space of braid matrices.
}

\noindent Observe that the Poisson bracket (3.2) does not depend on the times
 $u_1, u_2, \ldots, u_n$, although the monodromy map does.
\bigskip
\noindent{\bf Proof:}
{\it 1.} As a first step we explicitly write the restriction of   the bracket
(2.12) to the space of representations. By direct calculation  one
obtains
$$
\eqalignno{\{\tr (M_i M_k),&\tr (M_j M_l)\}=\{(M_i)^a_b(M_k)^b_a,
(M_j)^c_d(M_l)^d_c\}=\cr
&=(M_i)^a_b(M_j)^c_d\{(M_k)^b_a,(M_l)^d_c\}+
(M_i)^a_b(M_l)^d_c\{(M_k)^b_a,(M_j)^c_d\}+\cr
&+(M_k)^b_a(M_j)^c_d\{(M_i)^a_b,(M_l)^d_c\}+
(M_k)^b_a(M_l)^d_c\{(M_i)^a_b,(M_j)^c_d\}.&(3.3)
 \cr} 
$$
where we mean summation over repeated indices; using (3.1a), one can
rewrite the
left hand sides of (3.3) as  
$$
\{\tr (M_i M_k),\tr (M_j M_l)\}=\{n-4+s_{ik}^2,n-4+s_{jl}^2\}=
4s_{ik}s_{ij}\{s_{ik},s_{ij}\}.\eqno{(3.4)}
$$
Now one has to distinguish between three essentially different cases,  in
correspondence with the different order of the indices. 
\bigskip
\bigskip
\noindent $\bullet$  $i<k<j<l$ or $i<j<l<k$:
\bigskip
For $i<k<j<l$  all the addenda in  the right hand side of (3.3) 
involve a Poisson
bracket of the form
 (2.12) with correctly ordered indices. Here we write explicitly only the first
one: $$
i\pi\tr\biggl(M_iM_jM_lM_k+M_iM_kM_lM_j-M_iM_lM_jM_k-M_iM_kM_jM_l\biggr).
$$
The others have a similar form, and it is easy to see that they cancel
pairwise (the first with the second and the third with the fourth). 

The same happens when $i<j<l<k$, since the only 
difference  is a change of sign in the two last elements.
Hence it follows 
$$\eqalignno{\{\tr (M_i M_k),\tr (M_j M_l)\}&=0\qquad i<k<j<l&(3.5a)\cr
\{\tr (M_i M_k),\tr (M_j M_l)\}&=0\qquad i<j<l<k&(3.5b)\cr}
$$
 Using (3.1b) one
immediately obtains equations (3.2d/e) 
\bigskip
\bigskip
\noindent $\bullet$ $i<j<k<l$
\bigskip
Here the different order of the indices induces a change of sign in the second
addendum, which becomes  equal to the first. Equation (3.3) gives
$$
\eqalign{\{\tr (M_i M_k),\tr (M_j M_l)\}=&
2i\pi\tr\biggl(M_iM_jM_lM_k+M_iM_kM_lM_j-M_iM_lM_jM_k-M_iM_kM_jM_l\biggr)\cr
=&4i\pi s_{ik}s_{jl}(s_{ij}s_{kl}-s_{il}s_{kj}),\cr}
$$
where the last equality follows from Lemma 3.1.
Using eq.(3.4)
we obtain immediately eq. (3.2f)

\bigskip
\bigskip
\noindent $\bullet$ $i=j<k<l$ or $i<j=k<l$ or  $i<j<k=l$ 
\bigskip
If two 
indices coincide, for instance  $i=j<k<l$, the other two
cases are analogous,  we find
$$
\eqalignno{\{\tr (M_i M_k),&\tr (M_i M_l)\}=\{(M_i)^a_b(M_k)^b_a,
(M_i)^c_d(M_l)^d_c\}=\cr
&=(M_i)^a_b(M_i)^c_d\{(M_k)^b_a,(M_l)^d_c\}+
(M_i)^a_b(M_l)^d_c\{(M_k)^b_a,(M_i)^c_d\}+\cr
&+(M_k)^b_a(M_i)^c_d\{(M_i)^a_b,(M_l)^d_c\}+
(M_k)^b_a(M_l)^d_c\{(M_i)^a_b,(M_i)^c_d\}.\cr} 
$$
The first and the third addendum cancel, the last is zero (because $M_i^2=\id$), 
 and it
remains:
$$
\eqalign{\{\tr (M_i M_k),\tr (M_i
M_l)\}&=2i\pi\biggl((\tr(M_i^2M_lM_k)-\tr(M_iM_kM_iM_l)
\biggr)\cr
&=2i\pi[(n-4 +s_{kl}^2)-
(n-4+ s_{kl}^2+s_{ik}^2s_{il}^2-2s_{kl}s_{ik}s_{il})]\cr
&=2i\pi s_{ik}s_{il}(2
s_{kl}-s_{ik}s_{il}),\cr} 
$$
where the second equality follows from (3.1a) and (3.1b). Using (3.4) this leads
to (3.2a/b/c).

\bigskip
\bigskip
2. It follows from the commutativity of the diagram (1.7), where all the
arrows
are Poisson maps 

\bigskip
3. As we have said above, the eigenvalues of the monodromy matrices
are the Casimir functions for this Poisson
structure. Particularly, applying to $M_\infty$ we obtain, due to (2.3),
the needed statement. Practically it is more convenient to use
the coefficients
of the characteristic polynomial $det (S^{-1}S^T-\mu \id)$ as the basic
Casimirs.
\bigskip
4.  Recall [D], that the natural action of the braid group $B_n$ with $n$
strands on the space of Stokes matrices is generated by the following
transformations corresponding to the standard generators $\sigma_1$,
\dots, $\sigma_{n-1}$
$$
\sigma_i: ~ S\mapsto K_i S\, K_i
$$
where the matrix $K_i=K_i(S)$ has the form
$$
\eqalign{K_{jj}&=1,\quad j=1,\ldots,n;\quad j\not=i,i+1\cr
K_{ii}&=-s_{i i+1},\quad K_{i i+1}=K_{i+1 i}=1,\quad K_{i+1 i+1}=0.}
$$
Other matrix entries of $K_i$ vanish.
According to [D] this action describes the
structure of analytic continuation of the solutions of MPDE. Our Poisson
bracket is obviously invariant w.r.t. analytic continuation.

\hfill Q.E.D.
\bigskip
\bigskip
{\bf Example 1.} $n=3$. In this case the space of Stokes matrices has
dimension 3. Denoting $x=s_{12}$, $y=s_{13}$, $z=s_{23}$ we obtain, 
$$
\eqalign{
\{ x, y\} &={i\pi\over 2}(2z-xy)\cr
\{ y, z\} &= {i\pi\over 2}(2x-yz)\cr
\{ z,x\} &={i\pi\over 2}(2y-zx).\cr}
$$
Our Poisson bracket coincides, within the constant factor $-{i\pi\over 2}$, with
that of [D].
\bigskip
\bigskip
{\bf Example 2.} $n=4$. For convenience of the reader we write here down, 
omitting the constant factor ${i\pi\over 2}$, 
the Poisson bracket on the six-dimensional space of the Stokes
matrices of the form
$$
S=\left(\matrix{1&p&q&r\cr
0&1&x&y\cr
0&0&1&z\cr
0&0&0&1\cr}\right).
$$
\bigskip

$$
\eqalign{
\{p,q\}=(2x-pq)&\qquad\qquad\{x,y\}=(2z-xy)\cr
\{p,r\}=(2y-pr)&\qquad\qquad\{y,z\}=(2x-yz)\cr
\{q,r\}=(2z-qr)&\qquad\qquad\{z,x\}=(2y-zx)\cr}
$$

$$
\eqalignno{
\{x,p\}&=(2q-xp)\qquad\{q,x\}=(2p-qx)\qquad\{r,x\}=0\cr
\{y,p\}&=(2r-yp)\qquad\{q,y\}=2(pz-rx)\quad\ \{r,y\}=(2p-ry)\cr
\{p,z\}&=0\qquad\qquad\quad\ \{z,q\}=(2r-zq)\qquad\{r,z\}=(2q-rz)&(3.6)\cr}
$$
The Casimirs of this Poisson bracket are
$$
C_1=-4+p^2+q^2+r^2+
x^2+y^2+z^2-pqx-pry-qrz-
xyz+prxz
$$
and
$$
C_2=6-2(+p^2+q^2+r^2+
x^2+y^2+z^2)
+2(
-pqx
-pry-qrz-
xyz)-2(pqyz+qrxy)+
p^2r^2+q^2y^2+r^2x^2
$$

On the 4-dimensional level surfaces of the Casimirs the Poisson
bracket (3.6) induces a symplectic structure. These surfaces and the
symplectic structures on them are invariant w.r.t. the following action
of the braid group $B_4$:
$$
\eqalign{\sigma_1:(p,q,r,x,y,z)&\mapsto (-p,x-pq,y-pr,q,r,z)\cr
\sigma_2:(p,q,r,x,y,z)&\mapsto (q-px,p,r,-x,z-xy,y)\cr
\sigma_3:(p,q,r,x,y,z)&\mapsto (p,r-qz,q,y-xz,x,-z)\cr}
$$
\vfill\eject
\bigskip
\bigskip
\noindent {\scap Acknowledgments}
\bigskip
The author would like to thank Prof. B. Dubrovin for his guidance and for 
the careful reading of the manuscript.
\bigskip
\bigskip
\bigskip
\noindent{\scap References}
\bigskip
\noindent [A] M.Audin: {\it Lectures on Integrable Systems and Gauge theory}, In 
"Gauge theory and symplectic geometry", (Montreal, PQ, 1995), 1-48, NATO 
Adv.Sci.Inst.Ser.C Math.Phys.Sci. 488, Kluwer Acad. Publ. Dordrecht, (1997).
Preprint IRMA-1995/20 (1995).
\bigskip
\noindent [AM] A.Yu.Alekseev, A.Z.Malkin: {\it Symplectic Structure
of Moduli Space 
of Flat
Connection on a Riemann Surface}, Comm.Math.Phys. {\bf 169}, 99-119 (1995).
\bigskip
\noindent [ARS] M.J.Ablowitz, A.Ramani, H.Segur: {\it Nonlinear
Evolution Equations and Ordinary Differential Equations of
Painlev\'e Type}, Lett.Nuovo Cimento
{\bf 23}, 333-338, (1978).
\bigskip
\noindent\ \qquad  M.J.Ablowitz, A.Ramani, H.Segur: {\it A connection between 
nonlinear evolution equations and ordinary differential equations 
of Painlev\'e type. I and II }, J.Math.Phys. {\bf 76},
715-721 and 1006-1015, (1980).
\bigskip
\noindent [B] N.Bourbaki: {\it Groupes et Alg\`ebres de Lie}, Chap. 4,5 et
6,
Masson, Paris--New York--Barcelone--Milan--Mexico--Rio de Janeiro (1981).
\bigskip
\noindent [BJL1] W.Balser, W.B.Jurkat, D.A.Lutz: {\it Birkhoff invariants
and Stokes multipliers for meromorphic linear differential
equations }, J.Math.Anal.Appl
 {\bf 71}, 48-94 (1979).
\bigskip
\noindent [BJL] W.Balser, W.B.Jurkat, D.A.Lutz: {\it On the reduction
of connection problems for differential equations with an irregular
singular point to ones with only regular singularities}, SIAM J.Math.Anal.
 {\bf 12}, 691-721 (1981).
\bigskip
\noindent [D] B. Dubrovin:
{\it Geometry of 2D topological field theory}, in ``Integrable Systems and
Quantum Group", eds M.Francaviglia, S.Greco, Springer Lecture Notes
in Math. {\bf 1620}, 120-348, (1996). 
 \bigskip
\noindent [D1] B. Dubrovin:
{\it Painlev\'e transcendents in two--dimensional topological field theory},
Preprint SISSA 24/98/FM. To appear in Procedings of 1996 Carg\`ese
summer school `` Painlev\'e Transcendents: One Century Later". 
\bigskip
\noindent [FN1] H.Flaschka, A.C.Newell: {\it Monodromy--and 
Spectrum--Preserving--Deformations I}, Comm.Math.Phys. {\bf 76}, 65-116, (1980).
 \bigskip
\noindent [FN2] H.Flaschka, A.C.Newell: {\it The inverse monodromy transform
is a canonical transformation} in Math.Studies {\bf 61}, North Holland, 65-91,
 (1982).
\bigskip
\noindent [FR] V.V.Fock, A.A.Rosly:
{\it Poisson structure on moduli of flat connections on Riemann surfaces
and $r$--matrix}, Preprint ITEP-72-92, (1992).
\bigskip
\noindent [FT] L.D.Faddeev, L.A.Takhtajan: {\it Hamiltonian Methods in the 
Theory of Solitons } Springer Verlag, Berlin (1986).
\bigskip
\noindent [Fu] R.Fuchs: {\it Sur quelques \'equations diff\'erentielles
lin\'eaires du second ordre}, C.R.Acad.Sc.Paris, {\bf 141}, 555-558, (1905).
\bigskip
\noindent [G] R.Garnier: {\it Sur les \'equations diff\'erentielles
du troisi\`eme ordre dont l'integral g\'en\'erale est uniforme et sur
une classe d'\'equations nouvelles d'ordre superieur dont
l'integral g\'en\'erale a ses points critiques fix\'es} Ann.Sci.Ecole.Norm.Sup.
{\bf 29}, 1-126, (1912).
\bigskip
\noindent\ \quad R.Garnier: {\it Solution du probl\`eme de Riemann pour les 
sist\`emes diff\'erentiels dont l'integral g\'en\'erale est \`a  points 
critiques fix\'es} Ann.Sci.Ecole.Norm.Sup.
{\bf 43}, 177-307, (1926).
\bigskip
\noindent [Ha1] J.Harnad: {\it Dual Isomonodromic
Deformations and Moment Maps to Loop Algebras}, 
Comm.math.phys. 

\noindent  {\bf 166},337-365, (1994).
\bigskip
\noindent [Ha2] J.Harnad: {\it Quantum Isomonodromic
Deformations and the Knizhnik--Zamolodchikov Equations}, 

\noindent preprint
CRM-2890 (1994).
\bigskip
\noindent [Hi] N.Hitchin: {\it Frobenius manifolds}, in "Gauge Theory and 
symplectic geometry" (Montreal, PQ, 1995), 68-112, NATO 
Adv.Sci.Inst.Ser.C Math.Phys.Sci. 488, Kluwer Acad. Publ. Dordrecht, (1997).
\bigskip
\noindent [I] E.L.Ince: {\it Ordinary differential equations}, London--New 
York etc.,Longmans, Green and Co., 1972.

 \bigskip
\noindent [IN] A.R.Its, V.Yu.Novokshenov: {\it The isomonodromic deformation 
method in the 
theory of Painlev\'{e} equations}, Lecture Notes in Math. {\bf 1191}, 
Springer-Verlag, Berlin 1986.

 \bigskip
\noindent [JM] M.Jimbo, T.Miwa: {\it Monodromy preserving 
deformations of linear 
ordinary differential equations with rational coefficients. II.}
Physica {\bf 2D}, 407-448, (1981).
\bigskip
\noindent [JMU] M.Jimbo, T.Miwa, K.Ueno: {\it Monodromy  preserving
deformations of linear 
ordinary differential equations with rational coefficients. I.}
Physica {\bf 2D}, 306-352, (1981).
\bigskip
\noindent [KS] D.Korotkin, H.Samtleben:
{\it Quantization of Coset Space $\sigma$ --Models Coupled to 
Two--Dimensional Gravity}, Comm. Math. Phys. {\bf 190},411-457,  (1997).
\bigskip
\noindent [R] N.Reshetikhin: {\it The Knizhnik--Zamolodchikov System as a
deformation of the Isomonodromy Problem}, Lett. Math. Phys. {\bf 26}, 167-177,
 (1992).
\bigskip
\noindent [Sch] R.Sch\"afke:{\it \"Uber das globale
Verhalten der Normall\"osungen von $x^{'}(t)=
(B+t^{-1}A)x(t)$ und zweirer Arten von assoziierten Funktionen},
Math. Nachr. {\bf 121}, 123-145, (1985). 
\bigskip
\noindent [Schl] L.Schlesinger: {\it \"Uber eine Klasse von
Differentsialsystemen beliebliger Ordnung mit festen kritischer Punkten}
 J.f\"ur Math.{\bf 141}, 96-145, (1912).
\bigskip
\noindent [Si] Y.Sibuya: {\it Linear differential equations in the complex
domain: problems of analytic continuation}, AMS, Translations of
mathematical monographs, {\bf 82}, (1990).
\bigskip
\noindent [U] K.Ueno: {\it Monodromy preserving deformation of linear differential 
equations with irregular singular points} Proc.Japan Acad. ser.A {\bf 56},  97-102, 
(1980).
\bigskip
\noindent\ \quad K.Ueno: {\it Monodromy preserving deformation and its application to 
soliton theory} I and II Proc.Japan Acad. ser.A {\bf 56}, 103-108 
and 210-215,q (1980).
\bigskip
\noindent [ZF] V.E.Zakharov, L.D.Faddeev: {\it Korteweg--de Vries equation, a 
completely integrable 
Hamiltonian system}, Funct.Anal.Appl, {\bf 5}, 280-287,(1971).
\bigskip

\bye